# ON THE CONSTRUCTION OF EVEN ORDER MAGIC SQUARES

ABDULLAHI UMAR, Sultan Qaboos University, Al-Khod, P. C. 123, Muscat, OMAN

**Abstract**

The aim of this note is to introduce a fast new general method for the construction of double and single even order magic squares. The method for double even order magic squares is fairly straight-forward but some adjustment is necessary for single even order magic squares.

## 1. Introduction

For centuries magic squares have fascinated specialists and non-specialists alike, and still continue to do so. Basically, the study of magic squares is concerned with two problems: finding general methods for constructing all n by n magic squares or constructing all magic squares of a given class and enumerating how many distinct (excluding reflections and rotations) magic squares are there for a given order $n$. The latter seems to be the more difficult problem as very little progress has been made. In fact, there is a unique 3 by 3 magic square, exactly 880 (eight-hundred and eighty) 4 by 4 magic squares and for n > 4 it is still unknown how many magic squares are there, except for some subclasses.

On the problem of constructing magic squares, now there are several methods for constructing magic squares of any order (see references [1] – [4] for example). The fastest and general method for constructing odd order magic squares by consecutive numbering was found by the Middle Easterners. This method appears fully developed in Persian manuscripts of the 13[th] century [3]. However, the author is not aware of an analogue method (by consecutive numbering) for constructing even order magic squares and so we set out to fill this gap. The construction given here is an improved version of that in [5]. Notice also that it is necessary to classify even order magic squares into double and single even and the squares of the latter are more difficult to construct.

## 2. Basic facts and definitions

A *primitive magic square* (referred to as a magic square in what follows) of *order n* is a square consisting of the $n^2$ distinct numbers *1, 2, 3, ... , $n^2$* in $n^2$ subsquares such that the sum of each row, column and main diagonals adds up to the same total, $n(n^2 + 1)/2$.

A *double even order* magic square is one whose order is divisible by 4.

A *single even order* magic square is one whose order is divisible by 2 but not divisible by 4.

A pair of natural numbers *(a, b)* is $(\pm k)$-*complementary* if $a + b = n^2 + 1 \pm k$, where *n* is the order of the magic square. If $k = 0$, the pair is simply called *complementary* and, *a* and *b* are *complements* of each other. There are $n^2/2$ distinct complementary pairs in an even order magic square.

An even order magic square is *associated* if all its complementary pairs are equidistant from the centre of symmetry of the square, i. e., the point where all the lines of symmetry meet.

A magic square is *parallel* if all its complementary pairs are placed parallel to one another.

A magic square is *mixed* if it is neither parallel nor associated.

## 3. The construction of double even order magic squares

Since a magic square is a sort of symmetrical or balanced square we may therefore impose a simple symmetry on the even-odd number distribution, and as a result, a simple even-odd symmetrical picture is given by Figure 1.

| Odd Numbers | Even Numbers |
|---|---|
| Even Numbers | Odd Numbers |

**Fig. 1.** The general form of a double even square

Beginning with the pair of outermost columns, place the complementary pairs horizontally, from top to bottom in the following order:

$$(1, n^2), (2, n^2 – 1), (3, n^2 – 2), \ldots , (n, n^2 – n + 1),$$

observing the even-odd distribution symmetry in Figure 1. Next, consider the pair of next outermost columns and place the complementary pairs horizontally, from top to bottom in the following order:

$$(n + 1, n^2 – n), (n + 2, n^2 – n – 1), \ldots , (2n, n^2 – 2n + 1),$$

again, observing the even-odd distribution symmetry in Figure 1. Continue the process until you reach the pair of innermost columns, always observing the even-odd distribution symmetry in Figure 1. Interchanging the lower diagonal and bottom complementary pairs of the pair of innermost columns yields a parallel magic square. The construction of an 8 by 8 square is shown in Figure 2a with the parallel magic square shown in Figure 2b after the two interchanges have taken place as shown by * and ^. (36 was in the position of 32 and vice-versa before the interchange and similarly for 33 and 29.)

| 1 | 9  | 17 | 25 | 40 | 48 | 56 | 64 |
|---|----|----|----|----|----|----|----|
| 2 | 10 | 18 | 26 | 39 | 47 | 55 | 63 |
| 3 | 11 | 19 | 27 | 38 | 46 | 54 | 62 |
| 4 | 12 | 20 | 28 | 37 | 45 | 53 | 61 |
| 5 | 13 | 21 | 29 | 36 | 44 | 52 | 60 |
| 6 | 14 | 22 | 30 | 35 | 43 | 51 | 59 |
| 7 | 15 | 23 | 31 | 34 | 42 | 50 | 58 |
| 8 | 16 | 24 | 32 | 33 | 41 | 49 | 57 |

| 1  | 9  | 17 | 25  | 40  | 48 | 56 | 64 |
|----|----|----|-----|-----|----|----|----|
| 63 | 55 | 47 | 39  | 26  | 18 | 10 | 2  |
| 3  | 11 | 19 | 27  | 38  | 46 | 54 | 62 |
| 61 | 53 | 45 | 37  | 28  | 20 | 12 | 4  |
| 60 | 52 | 44 | 32* | 33^ | 21 | 13 | 5  |
| 6  | 14 | 22 | 30  | 35  | 43 | 51 | 59 |
| 58 | 50 | 42 | 34  | 31  | 23 | 15 | 7  |
| 8  | 16 | 24 | 36* | 29^ | 41 | 49 | 57 |

**Fig. 2**(a) Constructing an 8 by 8 magic square,  (b)  A complete 8 by 8 magic square

The numbers *1,2, ... , n* ( here *n = 8*) are placed consecutively down the first column, the numbers *n + 1, n + 2, ... ,2n* are placed consecutively, down the second column and continue the process until the numbers $(n/2 – 1)n + 1, (n/2 – 1)n + 2, \ldots, n^2/2$ are placed in the (n/2)-th column. Then starting at the last column place the numbers $n^2, n^2 – 1, \ldots , n^2/2$ (in descending order) in a similar fashion. Interchanging the odd and even numbers symmetrically within each row to adhere to Figure 1 and carrying out the interchanges shown by * and ^ gives the required square.

| R/C | 1 | 2 | 3 | … | m | m + 1 | … | n - 2 | n - 1 | n |
|---|---|---|---|---|---|---|---|---|---|---|
| 1 | 1 | n + 1 | 2n + 1 | … | p − n + 1 | p + n | … | 2p − 2n | 2p − n | 2p |
| 2 | 2p − 1 | 2p − n − 1 | 2p − 2n − 1 | … | p + n − 1 | p − n + 2 | … | 2n + 2 | n + 2 | 2 |
| 3 | 3 | n + 3 | 2n + 3 | … | p − n + 3 | p + n − 2 | … | 2p − 2n − 2 | 2p − n − 2 | 2p − 2 |
| 4 | 2p − 3 | 2p − n − 3 | 2p − 2n − 3 | … | p + n − 3 | p − n + 4 | … | 2n + 4 | n + 4 | 4 |
| . | . | . | . | … | . | . | … | . | . | . |
| . | . | . | . | … | . | . | … | . | . | . |
| . | . | . | . | … | . | . | … | . | . | . |
| m − 1 | m − 1 | 3m − 1 | 5m − 1 | … | p − m − 1 | p + m + 2 | … | 2p − 5m + 2 | 2p − 3m + 2 | 2p − m + 2 |
| m | 2p − m + 1 | 2p − 3m + 1 | 2p − 5m + 1 | … | p + m + 1 | p − m | … | 5m | 3m | m |
| m + 1 | 2p − m | 2p − 3m | 2p − 5m | … | p* | p + 1^ | … | 5m + 1 | 3m + 1 | m + 1 |
| m + 2 | m + 2 | 3m + 2 | 5m + 2 | … | p − m + 2 | p + m − 1 | … | 2p − 5m − 1 | 2p − 3m − 1 | 2p − m − 1 |
| . | . | . | . | … | . | . | … | . | . | . |
| . | . | . | . | … | . | . | … | . | . | . |
| . | . | . | . | … | . | . | … | . | . | . |
| n − 3 | 2p − n + 4 | 2p − 2n + 4 | 2p − 3n + 4 | … | p + 4 | p − 3 | … | 3n − 3 | 2n − 3 | n − 3 |
| n − 2 | n − 2 | 2n − 2 | 3n − 2 | … | p − 2 | p + 3 | … | 2p − 3n + 3 | 2p − 2n + 3 | 2p − n + 3 |
| n −1 | 2p − n + 2 | 2p − 2n + 2 | 2p − 3n + 2 | … | p + 2 | p − 1 | … | 3n − 1 | 2n − 1 | n − 1 |
| n | n | 2n | 3n | … | p + m* | p − m + 1^ | … | 2p − 3n + 1 | 2p − 2n + 1 | 2p − n + 1 |

**Fig. 3.** A 'complete' *n by n* double even parallel magic square.

The general result for an $n$ by $n$ magic square is shown in Figure 3, where $p = n^2/2$ and $m = n/2$. Each row of the $n$ by $n$ square consists of $n/2$ complementary pairs thus making a total of $n(n^2 + 1)/2$, the required sum. Each column consists of $n/4$ (+1)-complementary pairs and $n/4$ (−1)-complementary pairs thus making a total of $(n/4)(n^2 + 2) + (n/4)n^2 = n(n^2 + 1)/2$, the required sum. The sums of the entries in the two main diagonals are easily seen to be $(n/2)n^2$ and $(n/2)n^2 + n$. The interchange of the lower diagonal and bottom complementary pairs of the pair of innermost columns can be shown to increase the leading diagonal by $n/2$ and simultaneously decrease the other main diagonal by the same amount, thus, giving the required sum in each case. Hence the square obtained is a magic square. Moreover, by construction it is a parallel magic square.

4.  **The construction of single even order magic squares**

Since the row sum $n(n^2 + 1)/2$ is odd ($n/2$ is odd and $n^2 + 1$ is odd), Figure 1 must be modified. One possible modification is as shown in Figure 4. The two innermost rows are omitted together with the middle $2n$ numbers (of the numbers $1, 2, 3, \ldots , n^2$). The middle $2n$ numbers consists of the following $n$ complementary pairs:

$(n^2 − 2n)/2 + 1, (n^2 + 2n)/2; (n^2 − 2n)/2 + 2, (n^2 + 2n)/2 − 1; \ldots ; n^2/2, n^2/2 + 1.$

| Odd Numbers | Even Numbers |
|---|---|
| $(n/2)$-   th row ||
| $(n/2 + 1)$-   th row ||
| Even Numbers | Odd Numbers |

**Fig. 4.** The general form of a single even square

The general construction is in many ways similar to that of the double even case. The construction of a *10 by 10* magic square is shown in Figures 5a and 5b with the two centre rows omitted.

| 1 | 9 | 17 | 25 | 33 | 68 | 76 | 84 | 92 | 100 |
|---|---|----|----|----|----|----|----|----|-----|
| 2 | 10 | 18 | 26 | 34 | 67 | 75 | 83 | 91 | 99 |
| 3 | 11 | 19 | 27 | 35 | 66 | 74 | 82 | 90 | 98 |
| 4 | 12 | 20 | 28 | 36 | 65 | 73 | 81 | 89 | 97 |
|   |    |    |    |    |    |    |    |    |     |
|   |    |    |    |    |    |    |    |    |     |
| 5 | 13 | 21 | 29 | 37 | 64 | 72 | 80 | 88 | 96 |
| 6 | 14 | 22 | 30 | 38 | 63 | 71 | 79 | 87 | 95 |
| 7 | 15 | 23 | 31 | 39 | 62 | 70 | 78 | 86 | 94 |
| 8 | 16 | 24 | 32 | 40 | 61 | 69 | 77 | 85 | 93 |

**Fig. 5.** (a) Constructing a *10 by 10* magic square with the two centre rows omitted

| 1 | 9 | 17 | 25 | 33 | 68 | 76 | 84 | 92 | 100 |
|---|---|----|----|----|----|----|----|----|-----|
| 99 | 91 | 83 | 75 | 67 | 34 | 26 | 18 | 10 | 2 |
| 3 | 11 | 19 | 27 | 35 | 66 | 74 | 82 | 90 | 98 |
| 97 | 89 | 81 | 73 | 65 | 36 | 28 | 20 | 12 | 4 |
|   |    |    |    |    |    |    |    |    |     |
|   |    |    |    |    |    |    |    |    |     |
| 96 | 88 | 80 | 72 | 64 | 37 | 29 | 21 | 13 | 5 |
| 6 | 14 | 22 | 30 | 38 | 63 | 71 | 79 | 87 | 95 |
| 94 | 86 | 78 | 70 | 62 | 39 | 31 | 23 | 15 | 7 |
| 8 | 16 | 24 | 32 | 40 | 61 | 69 | 77 | 85 | 93 |

**Fig. 5.** (b) A *10 by 10* magic square with the two centre rows omitted

To complete the magic square, first we define a finite sequence of natural numbers involving the unused complementary pairs. Let $\{a_j\}$ be a finite sequence with common difference *d = 1* whose first and last terms are given by $(n^2 - 2n)/2 + 1$ and $n^2/2$, respectively. The *(n/2)*-th and *(n/2 + 1)*-th rows are filled as shown in Figure 6a and 6b. There are *n* pairs of cells to be filled from the left numbered 1 to *n*. The numbers $a_1, a_2, \ldots, a_m$ are placed in cells 1, 2, … , *m* of the *(n/2 + 1)*-th row, respectively; while the numbers $a_n, a_{n-1}, \ldots, a_{m+1}$ are placed in cells m + 1, m + 2, … , *n* of the *(n/2)*-th row, respectively. The blank cells consist of the complements of the numbers above or below them (within the two centre rows). Now, interchanging (column-wise) the numbers

$a_2, a_4, \ldots, a_{m-1}$ and $a_{m+3}, a_{m+5}, \ldots, a_{n-2}$

with their complements above or below can be shown to fix the sums of the two rows. Note that for *m = 3* (equivalently, *n = 6*) only one interchange is required. All the other rows consist of *n/2* complementary pairs each, thus giving the required sum. Each column consists of *(n – 2) /4* (+1)-complementary pairs, *(n – 2) /4* (–1)-com-plementary pairs and a single complementary pair (from the two centre rows) thus making a total of $n(n^2 + 1)/2$,

| Position | 1 | 2 | 3 | ... | m-2 | m-1 | m | m+1 | m+2 | m+3 | ... | n-3 | n-2 | n-1 | n |
|---|---|---|---|---|---|---|---|---|---|---|---|---|---|---|---|
| (n/2)th-row | | | | | | | | $a_n$ | $a_{n-1}$ | $a_{n-2}$ | | $a_{m+4}$ | $a_{m+3}$ | $a_{m+2}$ | $a_{m+1}$ |
| (n/2 +1)th-row | $a_1$ | $a_2$ | $a_3$ | ... | | $a_{m-1}$ | $a_m$ | | | | | | | | |

**Fig. 6.** (a) Filling the two centre rows of an *n by n* single even order magic square

| Position | 1 | 2 | 3 | ... | m-2 | m-1 | m | m+1 | m+2 | m+3 | ... | n-3 | n-2 | n-1 | n |
|---|---|---|---|---|---|---|---|---|---|---|---|---|---|---|---|
| (n/2)th-row | | $a_2$ | | | | $a_{m-1}$ | | $a_n$ | $a_{n-1}$ | | | $a_{m+4}$ | | $a_{m+2}$ | $a_{m+1}$ |
| (n/2 +1)th-row | $a_1$ | | $a_3$ | ... | $a_{m-2}$ | | $a_m$ | | | $a_{n-2}$ | | | $a_{m+3}$ | | |

**Fig. 6.** (b) The 'complete' two centre rows of an *n by n* single even order magic square

**(The blank cells consist of the complements of the numbers that appear at the top or bottom within the two centre rows)**

the required sum. The sums of the entries in each of the two main diagonals can be shown to be $n(n^2 + 1)/2 + 2$ and $n(n^2 + 1)/2 – 2$. Thus interchanging the complementary pairs
$(n^2/2 + 5n/2 – 3, n^2/2 – 5n/2 + 4)$ and $(n^2/2 + 5n/2 – 2, n^2/2 – 5n/2 + 3)$
(immediately below and immediately above the two centre rows) can also be shown to fix the sums of the two main diagonals to the required sum. Hence we obtain a mixed magic square. Figure 7 shows a complete mixed magic square of order 10, where the interchanged complementary pairs are indicated by * and ^.

| 1 | 9 | 17 | 25 | 33 | 68 | 76 | 84 | 92 | 100 |
|---|---|---|---|---|---|---|---|---|---|
| 99 | 91 | 83 | 75 | 67 | 34 | 26 | 18 | 10 | 2 |
| 3 | 11 | 19 | 27 | 35 | 66 | 74 | 82 | 90 | 98 |
| 97 | 89 | 81 | 72* | 65 | 36 | 29^ | 20 | 12 | 4 |
| 60 | 42 | 58 | 44 | 56 | 50 | 49 | 53 | 47 | 46 |
| 41 | 59 | 43 | 57 | 45 | 51 | 52 | 48 | 54 | 55 |
| 96 | 88 | 80 | 73* | 64 | 37 | 28^ | 21 | 13 | 5 |
| 6 | 14 | 22 | 30 | 38 | 63 | 71 | 79 | 87 | 95 |
| 94 | 86 | 78 | 70 | 62 | 39 | 31 | 23 | 15 | 7 |
| 8 | 16 | 24 | 32 | 40 | 61 | 69 | 77 | 85 | 93 |

**Fig. 7.** A complete *10 by 10* magic square

## 5. Concluding remarks.

The interesting thing about these methods of constructing even order magic squares is that they can be regarded as 'continuous' or 'consecutive' by analogy with the consecutive method for constructing odd order magic squares. For the double even case the description of the method is as follows:

'Place consecutively, the numbers 1, 2, … , $n$ from the first to the last row (one number in each row) alternating from left to right to left … in the pair of outermost columns. However, $n/2$ and $(n/2) + 1$ are to be placed in the same column, after which you resume the alternation from right to left to right …. Next place the numbers $n + 1, n + 2, …, 2n$ in exactly the same fashion as above in the pair of next outermost columns, noting that $3n/2$ and $(3n/2) + 1$ are to be placed in the same column. The process continues until the pair of innermost columns is reached and $p$ is placed in the bottom cell of the left-hand column of the innermost pair. Now place $p + 1$ adjacent to $p$ ( i. e., in the bottom cell of the right-hand column of the innermost pair), and continue consecutive placement but reversing the whole process. That is, place the numbers: $p + 1, p + 2, …, p + n$;

$p + n + 1, p + n + 2, …, p + 2n;$ … ; $2p – n + 1, 2p – n + 2, … , 2p$ alternately from right to left …, from the last to the first row and moving outwards from the pair of innermost columns to the pair of outermost columns. Note that at each half stage the left to right alternation pauses momentarily, then resumes after the half plus one stage. Finally, interchanging the lower diagonal and bottom complementary pairs of the pair of innermost columns completes the magic square '

e-mail: aumar@squ.edu.om